\documentclass[10pt, final]{amsart}
\usepackage{graphicx}
\usepackage{amssymb}

\newcommand{\C}{\mathbb{C}}
\newcommand{\R}{\mathbb{R}}
\newcommand{\Z}{\mathbb{Z}}
\newcommand{\abs}[1]{\lvert#1\rvert}
\newcommand{\norm}[1]{\lVert#1\rVert}
\newcommand{\paren}[1]{\left(#1\right)}
\newcommand{\bracket}[1]{\left[#1\right]}
\newcommand{\set}[1]{\left\{#1\right\}}

\usepackage{amsfonts}
\usepackage{amsmath}

\numberwithin{equation}{section}

\title{Numerical simulations of the energy- supercritical Nonlinear
Schr\"odinger equation}
\author{J. Colliander, G. Simpson, C. Sulem }

\begin{document}

\bibliographystyle{plain}

\begin{abstract}
We present numerical simulations of the defocusing nonlinear
Schr\"odinger (NLS) equation with an energy supercritical
nonlinearity.  These computations were motivated by recent works 
of Kenig-Merle and Kilip-Visan  who considered some energy supercritical wave equations and proved that if the solution is \emph{a priori} bounded in the critical Sobolev space (i.e.  the space whose homogeneous norm is invariant under the scaling leaving the equation invariant), then it exists for all time and scatters. 

 In this paper, we numerically investigate  the  boundedness of the  $H^2$-critical Sobolev 
norm for solutions of the NLS equation in dimension five with quintic nonlinearity.
 We find that for a  class of initial conditions, this norm  remains
 bounded, the solution exists for long time,  and  scatters.
\end{abstract}
\maketitle

\;

\section{Introduction}
This paper is a numerical investigation of wellposedness and scattering
properties  of solutions of the  defocusing
 nonlinear Schr\"odinger (NLS)  equation,
\begin{equation}
\label{eq:nls}
\mathrm{i} u_t + \Delta u - \abs{u}^{p-1} u =0, 
\quad u:  \R^{d}\times \R \to \C
\end{equation}
in the energy supercritical setting.
The notion of $H^s$-criticality is associated with the scaling transformation
$u(x,t) \to u_\lambda(x,t) = \lambda^{2/(p-1)} u(\lambda x,  \lambda^2 t)$
that leaves  the NLS equation invariant. We say that the problem is 
$H^s$-critical, if the homogeneous $H^s$-norm of the solution remains
unchanged under the above scaling.
We place ourselves in the energy  ($H^1$) supercritical regime by assuming
that the dimension and the nonlinearity are such that the critical  Sobolev exponent $s_c$  satisfies 
\begin{equation}
\label{eq:scrit}
s_c \equiv  \frac{d}{2} - \frac{2}{p-1} > 1.
\end{equation}
The two conserved quantities of the NLS equation, 
the mass
\begin{equation} \label{mass} 
M(u) = \int_{\R^d} |u(x,t)|^2 dx
\end{equation}
and the  energy
\begin{equation}
E(u) = \int_{\R^d} (|\nabla u(x,t)|^2  +\frac{1}{p+1} |u|^{p+1} )dx
\end{equation}
have indeed led to a special interest to mass-critical ($s_c=0$) and energy-critical ($s_c=1$)
problems.  

Significant progress on global well-posedness and scattering has been made since the pioneering work on scattering by Ginibre and Velo \cite{ginibre1979ocnb}, Lin and Strauss \cite{lin1978das}, and Strauss \cite{strauss81nst}.  In these early works, scattering was shown to hold for a range of energy subcritical configurations with finite variance.  Subsequent work \cite{ginibre1985ste} proved scattering in the energy space.  These results, and refinements, are collected in \cite{strauss1989nwe}, \cite{sulem1999nse}, \cite{cazenave2003sse}.

%{\bf{JC: Edit next paragraph; which equation(s) are referenced here?}}

For energy subcritical problems with $s_c <1$, global well-posedness
and scattering results have been obtained for $H^s$ data with $s$ near
1 and away from $s_c$. Bourgain \cite{bourgain1998ses} established
scattering in $H^s$, for radially symmetric data, for all $s \in
(11/3, 1)$.  This was refined by Colliander, Keel, Staffilani, Tao and
Takaoka \cite{colliander2003gea}, where the solution emerging from
general initial data  was shown to scatter for $s \in (4/5, 1)$. These results establish global well-posedness and scattering for some data which have \emph{infinite} energy.

Scattering results at critical regularity have been established
starting ten years ago.
A breakthrough result \cite{bourgain1999gwd} on the 3d energy critical
problem for radial data established a new strategy for proving
scattering for initial data of critical regularity. With other ideas
in \cite{colliander2004gwp}, Bourgain's induction-on-energy strategy resolved the case for
general data. The 4 dimensional defocusing energy critical case was
then established by Ryckman and Visan \cite{RV}  with the higher dimensional case
thereafter \cite{Visan}. Other important advances on the energy
critical problem appear in \cite{grillakis2000nse},
\cite{tao2005gwp}. 
Scattering was obtained \cite{KM_Inv} beneath the
natural threshold size for the focusing energy critical
problem. Moreover, ideas in \cite{KM_Inv, KM_Acta} simplified the
implementation of the energy critical strategy leading to advances on other model
equations. Building on these developments, scattering in the
$L^2$-critical problem for large
radial $L^2$ data was established in the defocusing
case \cite{killip3188cns} and for radial data under the ground state mass in the
focusing case \cite{killip2007mcn} .  
% In the defocusing energy critical case, Bourgain \cite{bourgain1999gwd}, Grillakis \cite{grillakis2000nse}, 
%Tao \cite{tao2005gwp} proved that for radial data in the energy
%space, there is scattering.    Colliander et
%al. \cite{colliander2004gwp} proved global well-posedness and
%scattering for the quintic problem in dimension $d=3$ for general
%data.  

Scattering is expected to hold for general large data of critical
regularity for \eqref{eq:nls}. To date, no such result has been
established except in the energy critical case. Under the assumption
of bounded critical norm, the defocusing $H^{1/2}$ case has been shown \cite{Kenig_Merle_Hhalf}
to scatter for large critical data.

%Extensions to higher dimensions were made by Visan \cite{visan2007dec}.

%Important progress has recently been made in the mass critical case for large radial data, removing the assumption of finite variance.  Killip, Tao and Visan \cite{killip3188cns} established scattering in $L^2$ for $d=2$ and Killip, Visan and Zhang \cite{killip2007mcn} extended this result to higher dimensions.  In the focusing case, if the initial mass is strictly smaller than that of the ground state, one also has scattering.

The hypothesis of a bounded critical Sobolev norm has recently been
considered in the energy supercritical regime $s_c > 1.$  Following a
recent work by Kenig and Merle \cite{kenig2008nrs} on the energy
supercritical wave{\footnote{The motivation for the numerical study
    described in this paper applies equally well to the
    nonlinear wave equation. As we had an existing code
    for simulating radially symmetric NLS , we chose to study NLS. In a future publication we shall
    examine the supercritical wave equation. Similar
    computations were performed by Strauss and Vazquez on nonlinear
    Klein-Gordon equations \cite{SV78}.}} equation, Killip and Visan \cite{killip2008esn}
considered some classes of  energy supercritical  NLS equations in dimension
$d\ge 5$ and proved that if the solution  is a priori  bounded in the critical Sobolev 
space $H^{s_c}$, it exits for all time and  scatters.

The purpose of our study  is to investigate
numerically
 the latter  assumption on the critical Sobolev norm. We are also
 motivated by the discussion of the ``theoretical possibility for
 computer assisted proofs of global well-posedness and scattering''
 appearing in \cite{B2000}.
We have performed our computations in the case 
$d=5$ with quintic nonlinearity ($p=5$). This is the `simplest'
case{\footnote{We observe similar phenomena in the $d=6, ~p=3$ case,
    which is also $H^2$ critical.}} when  the critical exponent is the smallest possible
 integer; $s_c=2$. In addition, we assume the initial
conditions are spherically symmetric 
to simplify the computations. 
Note that our choice of dimension and nonlinearity is not exactly
covered by Theorems 1.4 and 1.5 of \cite{killip2008esn}, although the authors claim that their result can be extended to any 
power law nonlinearity $|u|^{p-1} u$ with $p $ odd integer, and $s_c <p$.
%{\bf { C: I changed p to (p-1) to be consistent with our notations}}
 
The  main observation  of our numerical work is that,  for the various initial conditions we have
considered,  the critical norm $\dot H^2$ of the solution
remains bounded for all time and that the solution  scatters.
As time evolves and the solution reaches an asymptotic state, the
energy concentrates into the kinetic energy and the potential 
energy tends to zero. At the same time, the $\dot H^2$ norm  of the solution
stabilizes to some value.  This value may be much higher than that
of the initial conditions. 
For a more quantitative assessment of the solution in the frequency space,
we also calculate the norm of the solution in  the  Besov space $\dot{B}^2_{2,\infty}$. 
Let $P_j$ be the Fourier multiplier operator 
which projects the  Fourier transform of a function $f$  onto the annulus 
$2^j \leq |\xi| < 2^{j+1}$. 
The Besov $\dot B^s_{2,\infty} (\R^d)$ is equipped with the norm
\begin{equation} \label{besov-def}
 \sup_{j \in \Z}     2^{2j}    \| P_j f \|_{L^2_x} . 
\end{equation}
We observe that the $H^2$ density spreads out in Fourier space and the
$B^2_{2, \infty}$ norm shrinks to small values as time advances under
\eqref{eq:nls}. The emergent small Besov norm and Strichartz
refinements (such as those appearing in the work \cite{Planchon_Besov}
of Planchon and references therein) gives further evidence of
scattering for solutions of \eqref{eq:nls} and provides a smallness
mechanism for possibly implementing the computer assisted proof described in \cite{B2000}.

%{\bf{JC:COMMENTS ON BESOV AND BOURGAIN'PAPER}}

In Section 2 we describe our numerical schemes and in Section 3 we present our main results.

\section{Numerical Methods}
\subsection{Time and Space Discretization}
\label{sec:discrete}
We study \eqref{eq:nls} with $d=5$ and $p=5$ and radially symmetric data.  This configuration  is convenient because the $\dot H^2$  norm is bounded by the $L^2$ norm of $\Delta u$.   
%Since $\Delta u$ must be computed to integrate the equation, little additional work is required to approximate the norm.

To simulate  the problem, we first truncate the domain, restricting $r \in (0,R_{\max})$, with boundary conditions
\begin{equation}
u_r\mid_{r=0} = 0,\qquad
u(r = R_{\max})  = 0.
\end{equation} 
These are interpreted as a symmetry condition at the origin, and an infinite barrier at $r=R_{\max}$ so that $u = 0$ at all $r>R_{\max}$.  $R_{\max}$ must be sufficiently large to avoid boundary interaction.  The domain is discretized as
\begin{equation}
\label{eq:domain_disc}
0 = r_0 < r_1 < r_2<\ldots<r_j = j h <\ldots< r_N = L.
\end{equation}
Letting
$U_j(t) = u(r_j, t)$,
the fourth order spatial discretization is
\begin{equation}
\label{eq:nlsrad_disc}
\begin{split}
&\imath \dot{U}_j + \frac{-U_{j-2}+ 16 U_{j-1} - 30 U_j + 16 U_{j+1} - U_{j+2}}{12 h^2} \\
&\quad + \frac{4}{r_j}\frac{U_{j-2} - 8 U_{j-1} + 8 U_{j+1} - U_{j+2} }{12 h} = \abs{U_j}^4 U_j.
\end{split}
\end{equation}
%For $j = 0,1, N-1, N$, 
and the discretized boundary conditions are:
\begin{align}
&U_{-1} = U_1   \; ; \; 
U_{-2}  = U_2,\\
&U_{N+1}  = 0  \; ; \;
U_{N+2}  = 0.
\end{align}
 We  integrate in time using the classical  fourth order Runge-Kutta  scheme with $\Delta t = O( \Delta r^2)$ to ensure numerical stability.

\subsection{Fourier Transform}
We need to compute the Fourier transform of the solution to evaluate its norm in the Besov space.
The Fourier transform of a radial function in $\mathbb{R}^d$ at $k = \abs{\mathbf{k}}$ is
\begin{equation}
\label{eq:radft}
\hat{u}(k) = \frac{1}{k^\nu}\int_0^\infty u(r) J_\nu(k  r) r^{d/2} dr, \quad \nu = \frac{d-2}{2}.
\end{equation}
This formula is derived in many texts on the Fourier transform, including Stein and Weiss \cite{stein1971ifa}\footnote{Stein and Weiss used a different definition of the Fourier Transform.  Thus, their formula is slightly different.}.  Indeed,  for radially symetric functions,
\[
\hat{u}(\mathbf{k}) = \frac{1}{(2\pi)^{d/2}}\int_{\mathbb{R}^d} e^{-i \mathbf{k}\cdot \mathbf{x}} u(\mathbf{x}) d \mathbf{x} = \frac{1}{(2\pi)^{d/2}} \int_0^\infty \set{\int_{\partial B_0(1)}e^{-i k r \hat{\mathbf k}\cdot  \mathbf y} dS_{\mathbf y}} u(r) r^{d-1}dr.
\]
The kernel is
\[
\int_{\partial B_0(1)}e^{-i k r \hat{\mathbf k}\cdot \mathbf y} dS_{\mathbf y} = (2\pi)^{d/2} (kr)^{-(d-2)/2} J_{(d-2)/2}(kr),
\]
leading to expression  \ref{eq:radft}.  

Following Cree and Bones \cite{cree1993ane}, we approximate the integral as follows. For $k>0$, we use the trapezoidal rule at the grid points $r_j = j \Delta r$,  $\Delta r = R_{\max} / N$ and $j = 0, 1, \ldots, N$.  At $k=0$, we instead approximate the integral
\begin{equation}
\hat{u}(k=0) =\frac{1}{2^\nu \Gamma(1+\nu)}\int_0^\infty u(r)r^{\nu + d/2} dr.
\end{equation}
The Nyquist frequency, $K_{\max}$, is related to our discretization by the expression:
\[
K_{\max} R_{\max} = \frac{N}{2}.
\]
The Fourier transform is computed at the points $k_j = j \Delta k$, $\Delta k = 1/(2R_{\max})$ and $j = 0, 1, \ldots, N$.  As discussed in their article, Cree and Bones found this method to be robust, though it is slow.

\subsection{Besov Approximation}
\label{sec:besov_numerics}
We now approximate the Besov space norm
\begin{equation}
\norm{u}_{\dot{B}^2_{2,\infty}} =\sup_{j \in\mathbb{Z}} 2^{2j} \norm{\hat{u}}_{L^2([2^{j}, 2^{j+1}))}
\end{equation}
For this purpose, we first identify values of $j$  for which these integrals can be meaningfully computed by the trapezoidal rule.  Let
\begin{align}
j_{\min} &= \mathrm{ceil}(\log(4 k_1) / \log(2)) \\
j_{\max} & =\mathrm{floor}(\log(k_N) / \log(2))
\end{align}
We choose $j_{\min}$ to guarantee at least four grid points $< 2^{j_{\min}}$.    For $ j_{\min} \leq j \leq j_{\max}-1$, the integral
\begin{equation}
\label{eq:dyadicinterior}
 \norm{\hat{u}}_{L^2([2^{j}, 2^{j+1}))}^2= \int_{2^j}^{2^{j+1}} \abs{\hat{u}}^2 dk \approx q_j^2
\end{equation}
is computed by the trapezoidal rule.  We also compute 
\begin{equation}
\label{eq:dyadicsmall}
\norm{\hat{u}}_{L^2([0, 2^{j_{\min}}))}^2 =  \int_{0}^{2^{j_{\min}}} \abs{\hat{u}}^2 dk \approx q_{j_{\min}-1}^2
\end{equation}
and
\begin{equation}
\norm{\hat{u}}_{L^2([2^{j_{\max}}, K_{\max}])}^2 = \int_{2^j_{\max}}^{K_{\max}} \abs{\hat{u}}^2 dk\approx q_{j_{\max}}^2
\end{equation}
With these integrals in hand, 
\begin{equation}
\norm{u}_{\dot{B}^2_{2,\infty}} \approx \max_{j_{\min}-1 \leq j \leq j_{\max}} \set{ 2^{2j} q_j}
\end{equation}

\subsection{Error of Numerical Scheme}

\begin{figure}
\begin{center}
\includegraphics[width=4in]{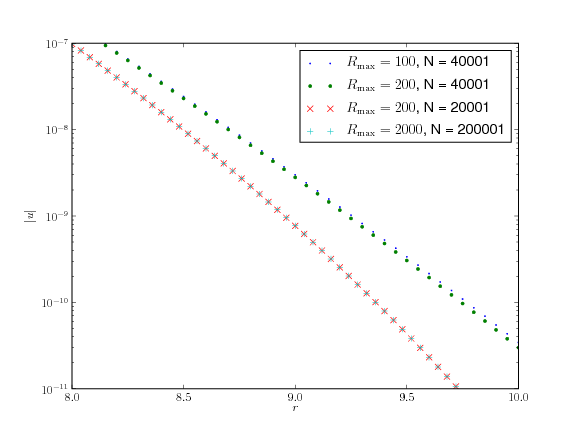}
\caption{Variation in the tails of $\abs{u}$ at $t=.04$ for Gaussian data $u_0= 10e^{-r^2}$.}
\label{fig:tailerror}
\end{center}
\end{figure}

We have tested our scheme by varying both the domain size and the grid resolution.  In Table \ref{table:amperror}, we show two metrics of our simulations, the value of $\abs{u}$ at the origin, and the maximum of $\abs{u}$ at a fixed time.  We see consistency amongst the simulations for the different parameters.  Examining the tails of $\abs{u}$ in Figure \ref{fig:tailerror}, we get a qualitative assessment of how these parameters influence the simulation.  So long as we have not reached the edge of the domain, which none of these simulations have, $R_{\max}$ matters little.  In the better resolved simulations, $\abs{u}$ has propogated farther too the right.  This is expected since greater resolution resolves higher wave numbers.  As will be argued in the next section, $u$ is scattering, and thus obeys the linear equation, where higher frequencies propagate with greater speed.

\begin{table}
\begin{center}
\begin{tabular}{l | l | l| l}
No. Points & $R_{\max}$ & $\abs{u}(r=0)$ & $\max_{r\in[0,R_{\max}]} \abs{u}$  \\
\hline
10000 & 100 & 0.7126025579 & 2.665689301 \\
20000 & 100 & 0.712561663 & 2.665668567 \\
40000 & 100 & 0.7125588732 & 2.665667313 \\
20000 & 200 & 0.7126025586 & 2.665689301 \\
40000 & 200 & 0.7125616583 & 2.665668567 \\
200000 & 2000 & 0.7126025579 & 2.665689301 
\end{tabular}
\caption{Convergence of $\abs{u}(r=0)$ and $\max_r \abs{u}$ for Gaussian data $u_0 = 10e^{-r^2}$ at $t = .02$.}
\label{table:amperror}
\end{center}
\end{table}

The discretization scheme is neither mass nor energy conservative. A calculation of the relative error on these conserved quantities is another measure of accuracy.   Table \ref{table:conlaw} shows that
the spatial resolution for our various initial conditions and  the relative error for  the two invariants.  As expected, the simulations with better spatial resolution, hence better temporal resolution, better conserve the invariants.  These invariants were computed by Simpson's method on the interval $[0,R_{\max}]$ with the discrete densities
\[
\bracket{(\Re U_i )^2 +(\Im U_i)^2}r_i^{4}
\]
and
\[
\bracket{-(\Re U_i \Delta_{\textrm{disc.}} \Re U_i) - (\Im U_i \Delta_{\textrm{disc.}} \Im U_i)  + \frac{2}{p+1} \paren{(\Re U_i )^2 +(\Im U_i)^2}^{6} }r_i^{4}.
\]
$\Delta_{\textrm{disc.}}$ is the discrete Laplacian from \eqref{eq:nlsrad_disc}.

\begin{table}
\begin{center}
\begin{tabular}{l | l | l | l | l | l}
Initial Condition & No. Points & $R_{\max}$ & $T_{\max}$ &  $\max_t\abs{\%\Delta\mathrm{Mass}}$ & $\max_t\abs{\%\Delta\mathrm{Energy}}$ \\
\hline
Gaussian, $u_0 = 10 e^{-r^2}$ & 40001 & 100 & 0.04 & 6.55478e-09 & 2.67077e-08 \\
Gaussian, $u_0 = 10 e^{-r^2}$ & 200001 & 2000 & 3.2 & 1.97657e-06& 7.94582e-06 \\
Ring, $u_0 = 8 r^2e^{-r^2}$ & 32001 & 100 & 0.2 &2.33731e-10&2.19067e-09 \\
Ring, $u_0 = 8 r^2e^{-r^2}$ & 120001 & 2400 & 9.0 & 9.61835e-07 &8.70909e-05 \\
Osc. Gaussian, $u_0 = 4 e^{-10\mathrm{i} r^2}e^{-r^2}$ & 40001 & 100 & 0.1 & 1.35194e-08 & 1.8279e-08 \\
Osc. Gaussian, $u_0 = 4 e^{-10\mathrm{i} r^2}e^{-r^2}$ & 200001 & 1000 & 1.0 & 2.15647e-07 & 2.91404e-07 
\end{tabular}
\caption{Error in the conserved quantities for each simulation.}
\label{table:conlaw}
\end{center}
\end{table}

\begin{figure}
\begin{center}
\includegraphics[width=4in]{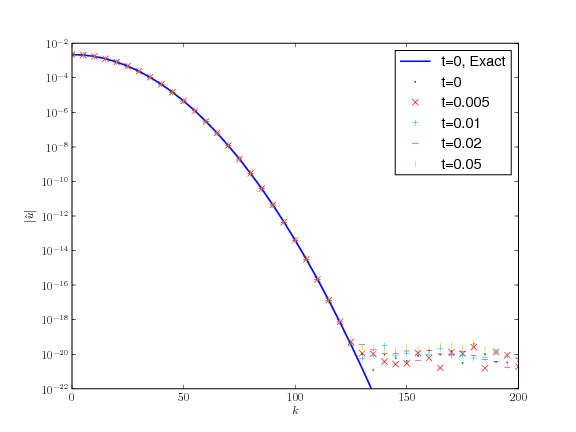}
\caption{Fourier transform at various times of a linear Schr\"odinger equation simulation with initial condition $u_0 = 4 e^{-10 \mathrm{i} r^2} e^{-r^2}$}
\label{fig:fouriercheck}
\end{center}
\end{figure}

We also verify our time stepping and Fourier approximation by simulating a \emph{linear} problem, computing the approximate Fourier transform, and observing that it does not change in time; see  Figure \ref{fig:fouriercheck}.

\section{Numerical Observations}
In this section we present and discuss our numerical simulations for energy supercritical defocusing NLS, \eqref{eq:nls}.  Throughout, our initial conditions are radially symmetric, simplifying the computations.  We speculate that the dynamics persist for general data and for other energy supercritical configurations.

\subsection{Initial Conditions}

We consider several families of initial conditions.  These are:
\begin{description}
\item[Gaussians] 
\begin{equation}
\label{eq:gaussiandata}
u_0(r) = A e^{-r^2}
\end{equation}
Under the linear flow, these are well known to spread and decay in $L^\infty$.  Our simulations show the nonlinear distortion in the shape.
\item[Rings]
\begin{equation}
\label{eq:ringdata}
u_0(r) = A r^2 e^{-r^2}
\end{equation}
\item[Oscillatory Gaussians]
\begin{equation}
u_0(r) = A e^{-\alpha \mathrm{i} r^2} e^{-r^2}
\end{equation}
Under the linear flow with $\alpha >0$, such an initial condition will initially focus towards the origin, then relax and decay.  Our simulations show that the nonlinearity arrests this focusing.
\end{description}

In all cases, the amplitude $A$ is taken ``large enough" so that the nonlinear effects in \eqref{eq:nls} will, at least initially, be strong and we will be outside the small data regime where scattering is known to occur \cite{sulem1999nse}.  We present the following cases: Gaussian data with $A=10$, $u_0 = 10 e^{-r^2}$; Ring data with $A=8$, $u_0 = 8 r^2 e^{-r^2}$; Oscillating Gaussian data with $A=4$ and $\alpha = 10$, $u_0 = 4 e^{-10 \mathrm{i} r^2} e^{-r^2}$.

\subsection{Results}

In all simulations, we find that after a transient period, the solution remains smooth and monotonically decays in amplitude.  This is evidence of global well-posedness.  The $L^\infty$ smallness of the solution is a first indication of scattering, as the nonlinearity becomes a perturbation of the linear equation.  After this transient period, the $L^6$-norm also decays monotonically.  Since the energy invariant is conserved, the potential energy is absorbed by the kinetic term, as is expected in scattering.  We also study the evolution of the $L^{14}_x$-norm, which is important because, by Strichartz,
\[
\norm{e^{i\Delta t} u_0}_{L^{14}_{t,x}}\leq C \norm{u_0}_{\dot H^2}.
\]
If the flow is to become asymptotically linear, we would expect $L^{14}_x$ to decay $\propto t^{-15/7}$, the theoretical rate of the linear flow. Indeed, when the simulation is run for a sufficiently long time, this is observed.  Finally, there is the critical norm, $\dot H^2$.  For numerical convenience we track $\norm{\Delta u}_{L^2}$, which controls $\dot H^2$.  Finally, our simulations indicate that $\norm{\Delta u}_{L^2}$  saturates to a finite value.  Thus, we have evidence that the scale invariant norm is globally bounded in time.  This was the fundamental a assumption in  \cite{kenig2008nrs, killip2008esn}.

% and the critical norm saturates.  This migration of energy into the kinetic state is evidence that the flow is becoming asymptotically linear; there is scattering.  We also observe that the $L^{14}_x$ norm, whose importance we will explain, decays  $\propto t^{-15/7}$.  This is the theoretical decay rate of the linear flow \cite{cazenave2003sse}\footnote{This linear result assumes the data has finite variance, which we appear to have.}.

\begin{figure}
\begin{center}
$\begin{array}{cc}
\includegraphics[width=3.5in]{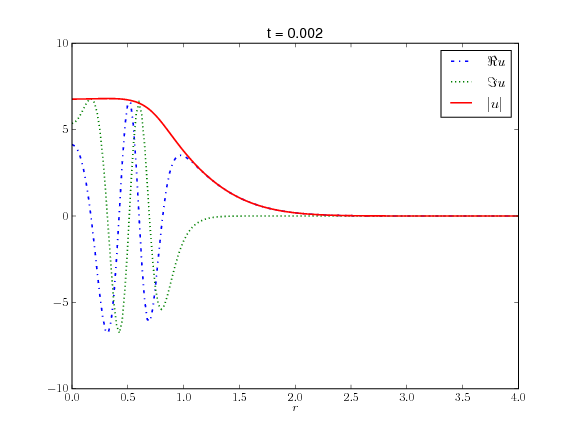}  &  
\includegraphics[width=3.5in]{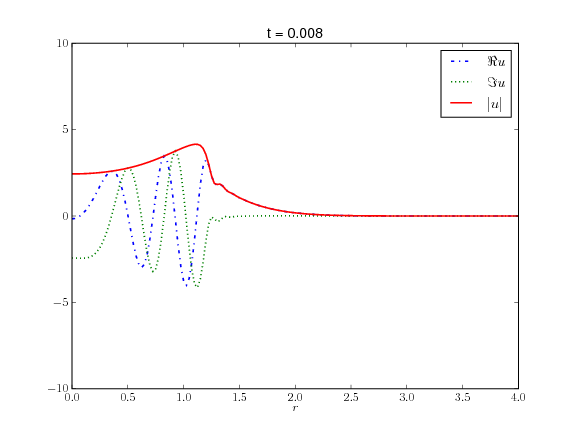}\\
\includegraphics[width=3.5in]{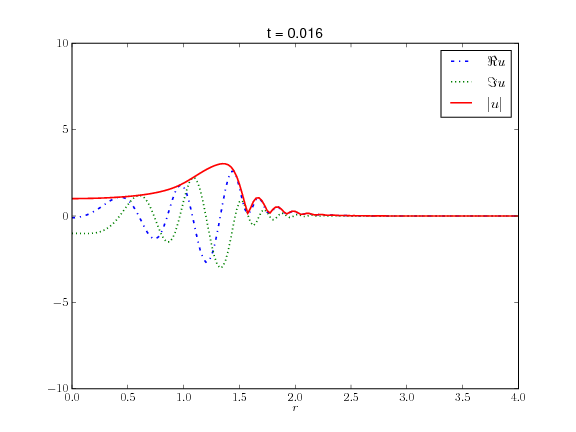}  &  
\includegraphics[width=3.5in]{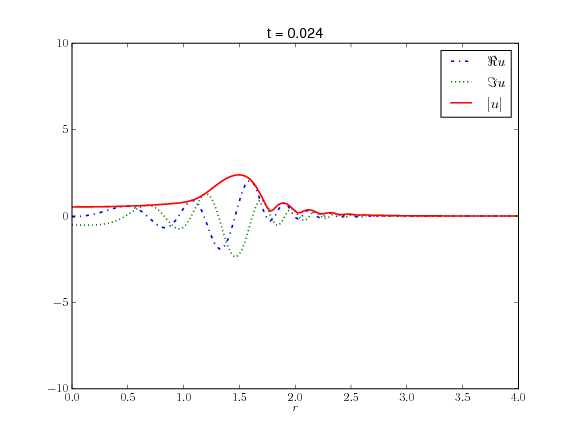}\\
\includegraphics[width=3.5in]{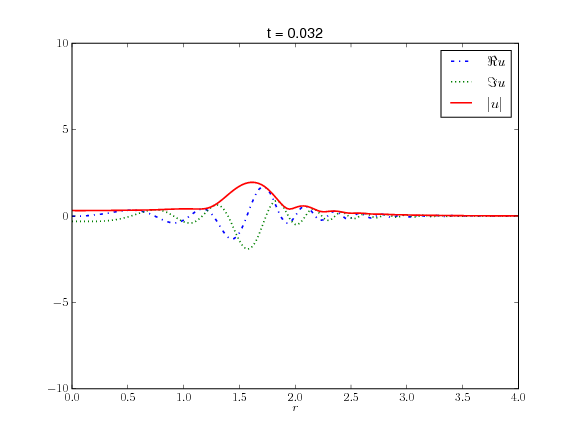} &
\includegraphics[width=3.5in]{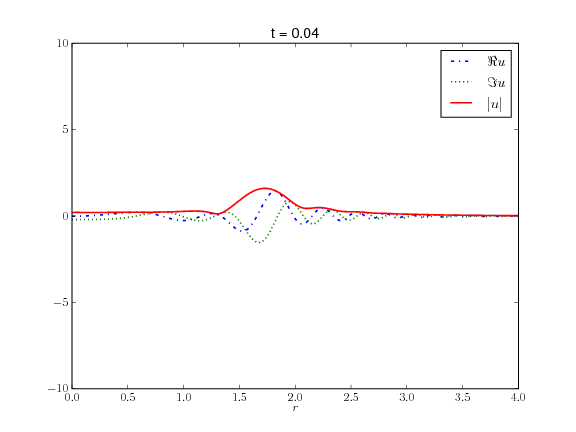}
\end{array}$
\caption{Evolution of $u_0 = 10 e^{-r^2}$.  Computed on the domain $[0, 100]$ with $40000+1$ points.}
\label{fig:gauss10profiles}
\end{center}
\end{figure}

\begin{figure}
$\begin{array}{cc}
\includegraphics[width=3.5in]{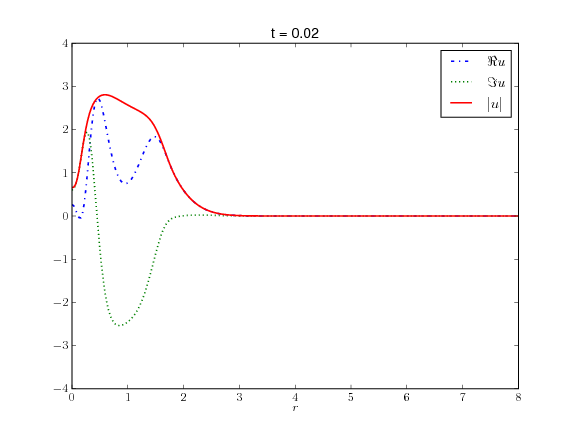}  &  
\includegraphics[width=3.5in]{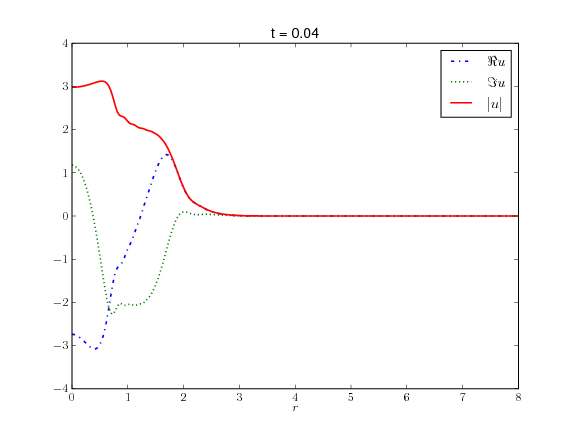}\\
\includegraphics[width=3.5in]{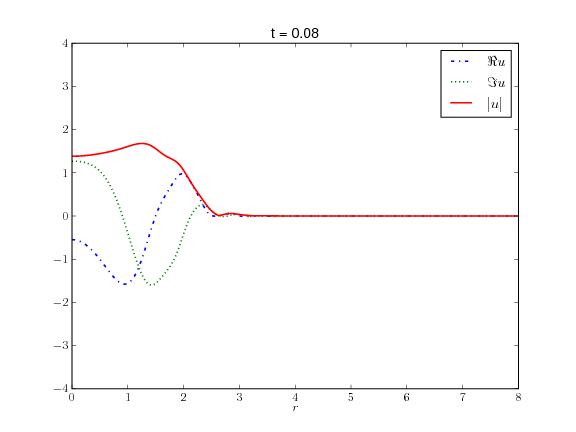}  &  
\includegraphics[width=3.5in]{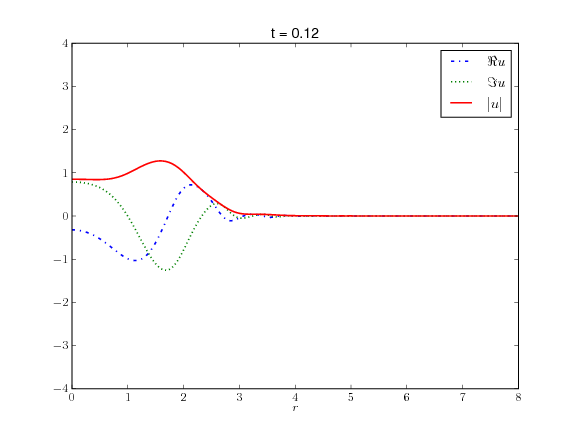}\\
\includegraphics[width=3.5in]{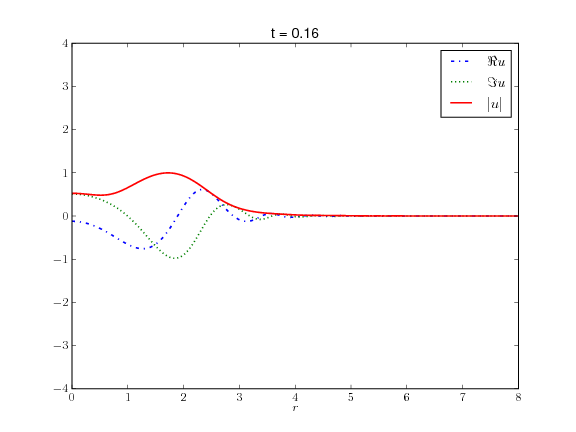} &
\includegraphics[width=3.5in]{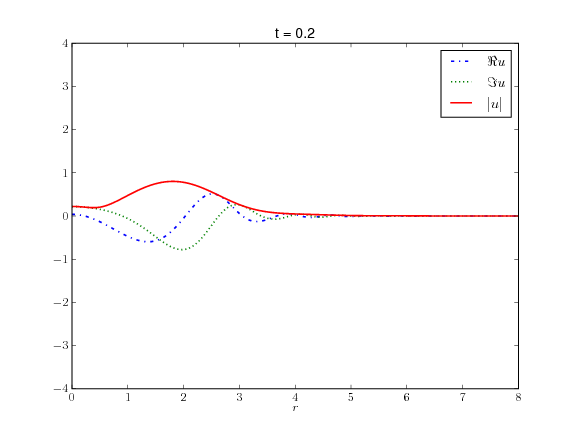}
\end{array}$
\caption{Evolution of $u_0 = 8 r^2 e^{-r^2}$.  Computed on the domain $[0, 100]$ with $32000+1$ points.}
\label{fig:ring7profiles}
\end{figure}

\begin{figure}
$\begin{array}{cc}
\includegraphics[width=3.5in]{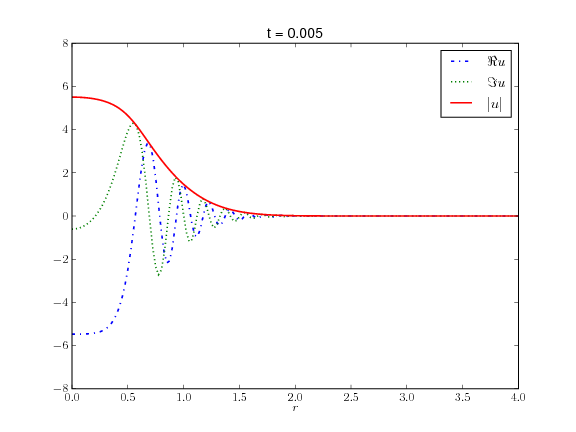}  &  
\includegraphics[width=3.5in]{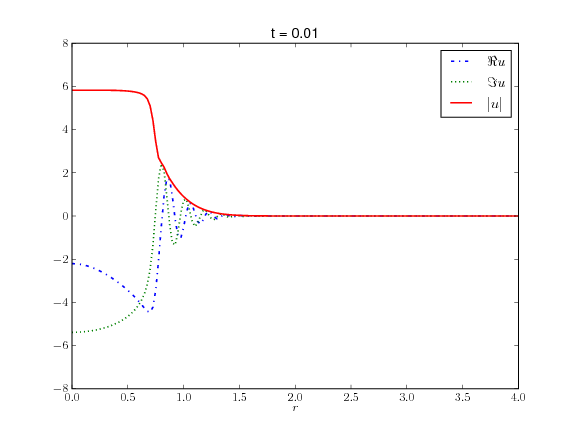}\\
\includegraphics[width=3.5in]{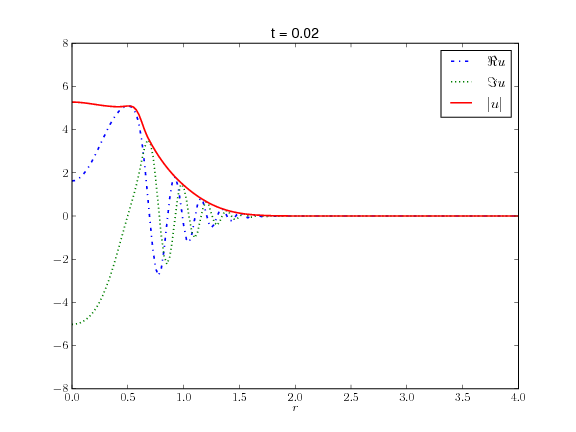}  &  
\includegraphics[width=3.5in]{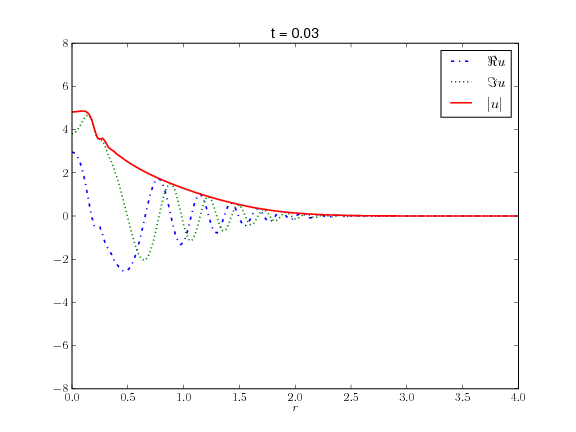}\\
\includegraphics[width=3.5in]{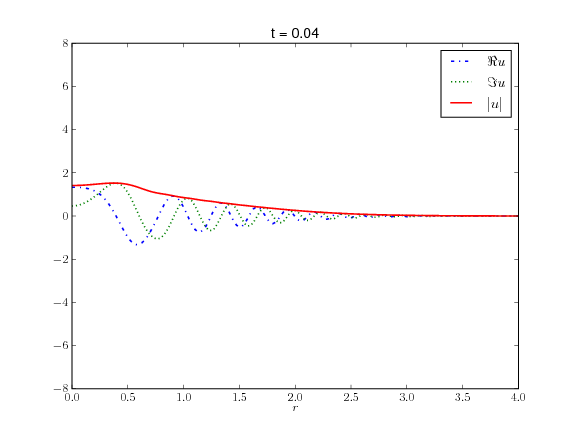} &
\includegraphics[width=3.5in]{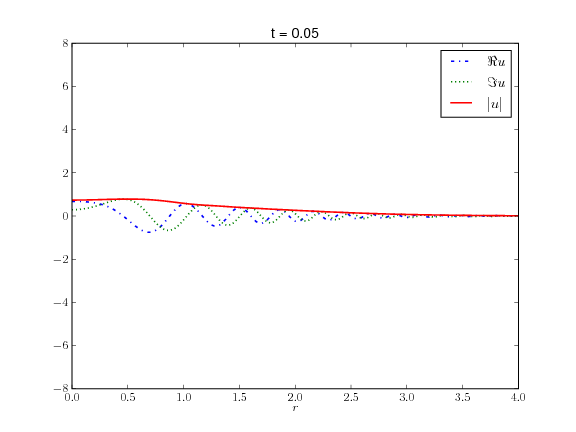}
\end{array}$
\caption{Evolution of $u_0 = 4 e^{-10\mathrm{i}r^2} e^{-r^2}$.  Computed on the domain $[0, 100]$ with $40000+1$ points.}
\label{fig:ring5profiles}
\end{figure}

Let us examine the profiles from our simulations.  The evolution of the Gaussian data, $u_0 = 10 e^{-r^2}$ is plotted in Figure \ref{fig:gauss10profiles}.  The shape is distorted, but $\abs{u}$ is monotonically decreasing in time.  We can also see that oscillations in the shape spread. The ring data, $u_0 = 8 r^2 e^{-r^2}$, has more complex transient dynamics.  As shown in Figure \ref{fig:ring7profiles}, there is initially a focusing of $u$ towards the origin.  This subsequently relaxes, and $u$ spreads and decays in amplitude, much like the Gaussian data.  The oscillatory gaussian, $u_0 = 4 e^{-10 \mathrm{i} r^2} e^{-r^2}$, is similar.  It initially focuses, subsequently relaxes, and appears to scatter as in Figure \ref{fig:ring5profiles}.  

\begin{figure}
$\begin{array}{cc}
\includegraphics[width=3.5in]{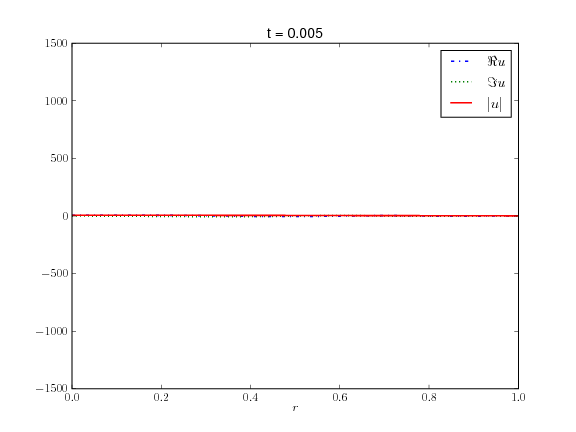}  &  
\includegraphics[width=3.5in]{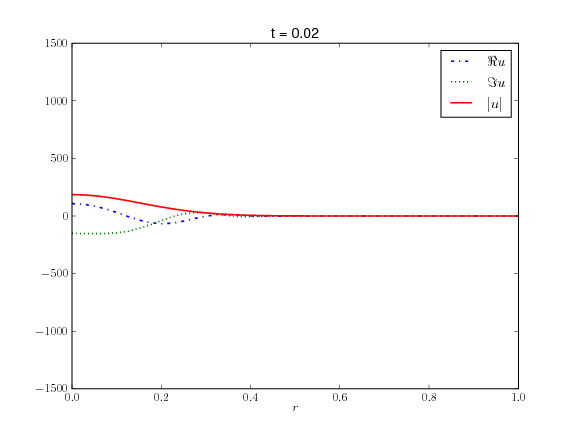}\\
\includegraphics[width=3.5in]{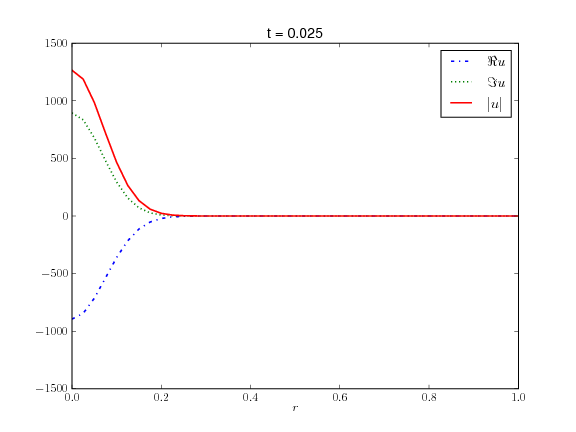}  &  
\includegraphics[width=3.5in]{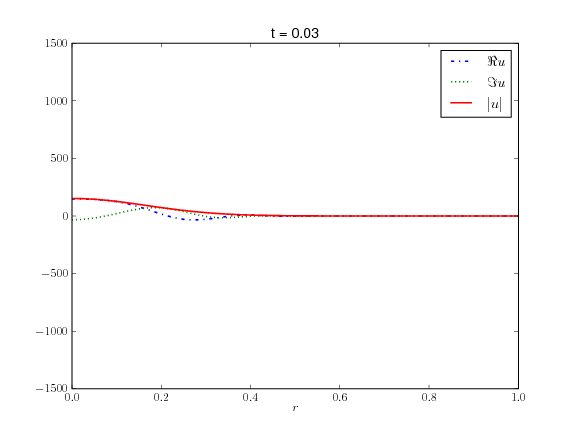}\\
\includegraphics[width=3.5in]{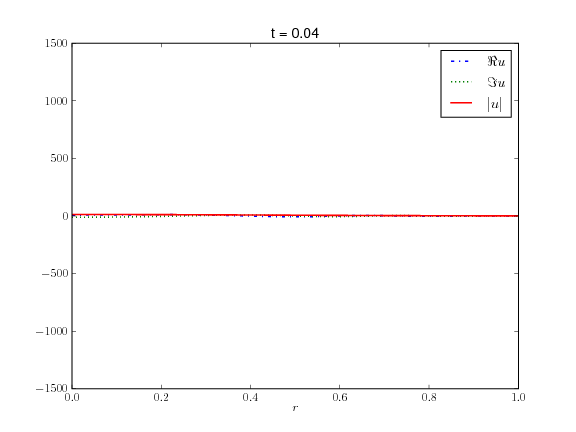} &
\includegraphics[width=3.5in]{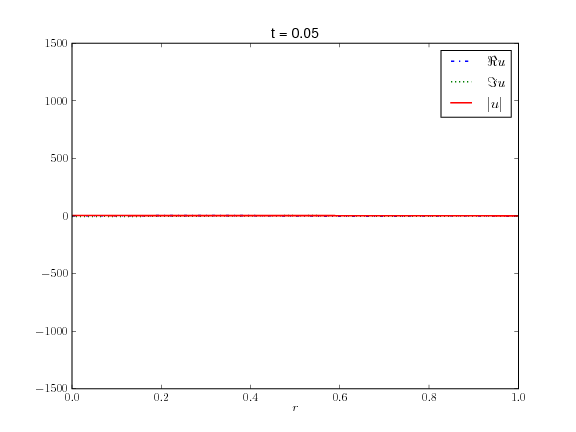} 
\end{array}$
\caption{Evolution of $u_0 = 4 e^{-10\mathrm{i}r^2} e^{-r^2}$ under the linear flow.  Computed on the domain $[0, 100]$ with $40000+1$ points.  }
\label{fig:ring5linearprofiles}
\end{figure}

For comparison, we simulate the linear Schr\"odinger equation with data $u_0 = 4 e^{-10 \mathrm{i} r^2} e^{-r^2}$ and the same discretization as in the nonlinear problem.  Figure \ref{fig:ring5linearprofiles} shows much stronger focusing towards the origin than seen in  Figure \ref{fig:ring5profiles}.  The nonlinearity is arresting the rush towards the origin and keeps the amplitude orders of magnitude smaller.

%have a similar dynamics as in the case of the Gaussian data.  We interpret this spreading and decrease in amplitude as scattering.   

%Consider the gaussian data $u_0 = 10 e^{-r^2}$.  Examining the short time evolution in Figure \ref{fig:gauss10profiles}, we see a monotonic, in time, decrease in the amplitude of $u$.  We also see the distortion associated with the nonlinearity; this is not a spreading and shrinking Gaussian.  

%For more interesting transient behavior, let us consider the ring data, $u_0 = 8 r^2 e^{-r^2}$.  As shown in Figure \ref{fig:ring7profiles}, there is initially a focusing of $u$ towards the origin.  This subsequently relaxes, and have a similar dynamics as in the case of the Gaussian data.  We interpret this spreading and decrease in amplitude as scattering.   

%The exact \emph{linear} solution for such data is:
%\begin{equation}
%u(r,t) = \frac{4}{(1 - (40 - 4\mathrm{i})t)^{5/2}}e^{-\frac{10 - \mathrm{i}}{-\mathrm{i} + (4 + 40 \mathrm{i})t}r^2}
%\end{equation}
%It is interesting that under the linear flow, this has a peak amplitude of $1000--1500$(I WILL WORK THIS OUT EXACTLY), far in excess of what is seen in the nonlinear simulations.

For a more quantitative assessment of scattering, we examine the aforementioned integrated quantities.  In Figure \ref{fig:gauss10metrics} (a) and (b), we plot $\norm{\Delta u}_{L^2}$.  Figure \ref{fig:gauss10metrics} (a) shows the rapid growth of this norm to $\sim 1200$, orders of magnitude larger than the initial value.   Figure \ref{fig:gauss10metrics} (b), computed for a longer time, suggests that it has saturated at this value.  As expected, the  potential energy, $\int \abs{u}^6$, vanishes.  We see this in in Figures \ref{fig:gauss10metrics} (c) and (d), where we have plotted $\norm{u}_{L^6}$ from the same two simulations.  Lastly, in Figure \ref{fig:gauss10metrics} (f) the $L^{14}$ space norm, after a sufficient time, begins to decay as $ \propto t^{-15/7}$.  

The ring data and the oscillatory gaussian are, asymptotically, very similar.  The same plots appear in Figures \ref{fig:ring7metrics} and \ref{fig:ring5metrics}.  We see rapid saturation of $\dot H^2$, the decay of the potential energy, and the asymptotically linear decay of the $L^{14}_x$ norm.  One notable difference is that in the oscillating gaussian simulation, the initial focusing causes a \emph{decrease} in $\dot H^2$ and an \emph{increase} in the potential.  Furthermore, the saturated value of $\norm{\Delta u}_{L^2}$ is not appreciably larger than in its initial value.  These simulations suggest there are at least two different time scales of interest.  Saturation of $\dot H^2$ happen very rapidly.  In contrast, the expected asymptotic decay of $L^{14}_x$ sets in at a much later time.

%
%This saturation of $\dot H^2$ is indicative of the flow becoming asymptotically linear.  As it becomes linear we would expect the energy to become entirely kinetic and the potential, $\int \abs{u}^6$, to vanish.  Indeed, We see exactly that in Figures \ref{fig:gauss10metrics} (c) and (d), where we have plotted $\norm{u}_{L^6}$ from the same two simulations.

%This is consistent with the conjecture that we are entering the asymptotically linear flow.  This norm is particularly relevant because through Strichartz,
%\[
%\norm{e^{i\Delta t} \abs{u}^4 u}_{\dot H^2}\leq C \norm{u}_{L^{14}_{t,x}}.
%\]
%This naturally appears when studying the evolution of $\dot H^2$ with Duhamel's principle.

\subsection{Fourier and Besov}
\begin{figure}
$\begin{array}{cc}
\includegraphics[width=3.5in]{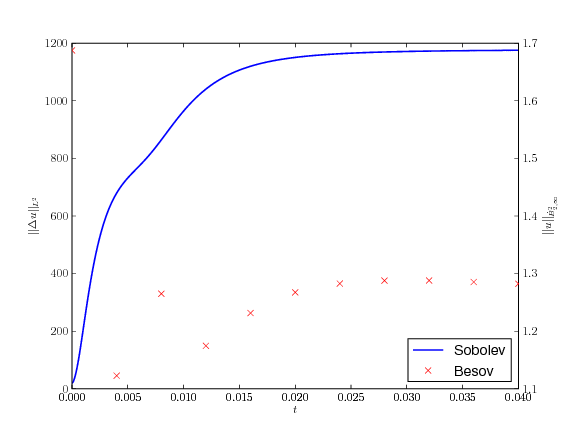}  &
\includegraphics[width=3.5in]{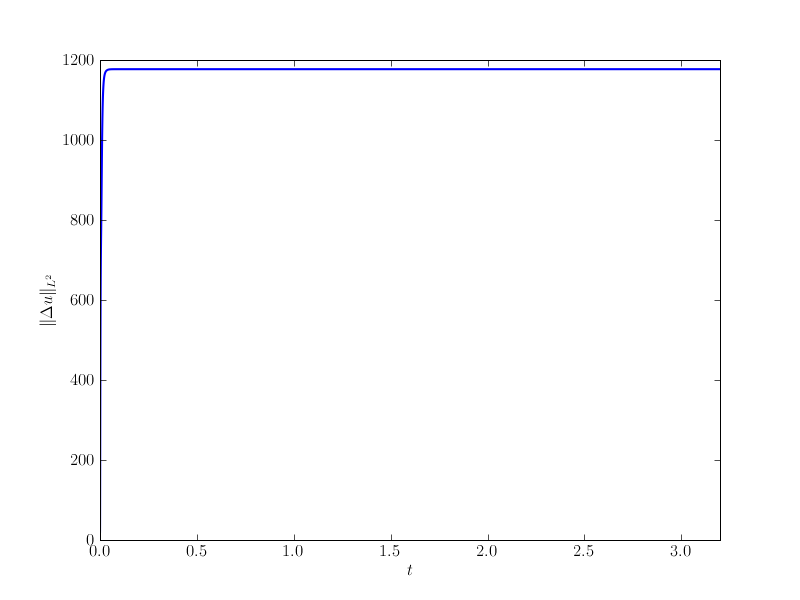}\\
\mathrm{(a)} & \mathrm{(b)}\\
\includegraphics[width=3.5in]{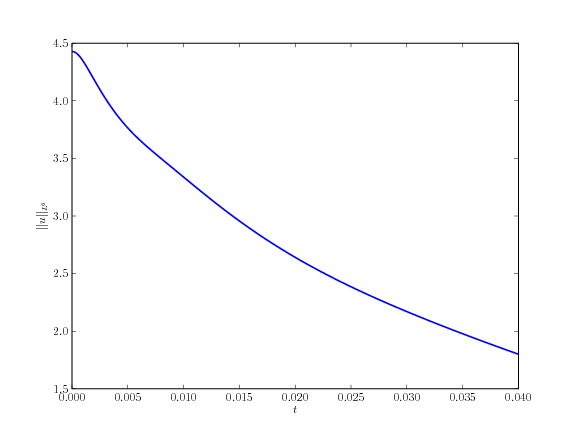}  &
\includegraphics[width=3.5in]{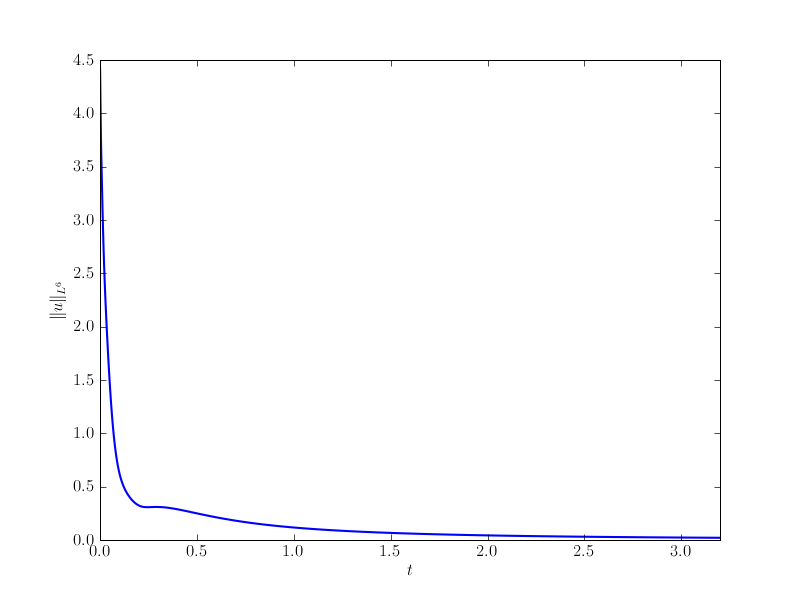}\\
\mathrm{(c)} & \mathrm{(d)}\\
\includegraphics[width=3.5in]{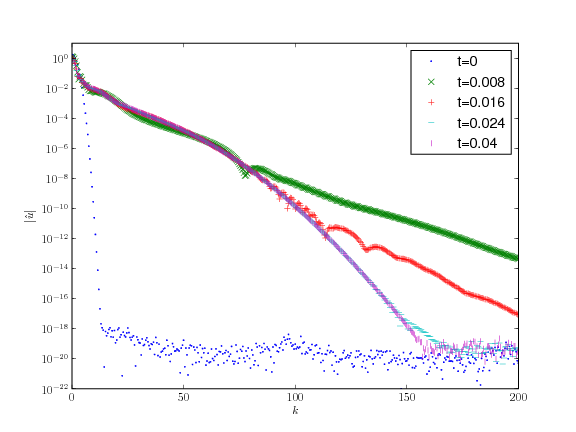}  &  \includegraphics[width=3.5in]{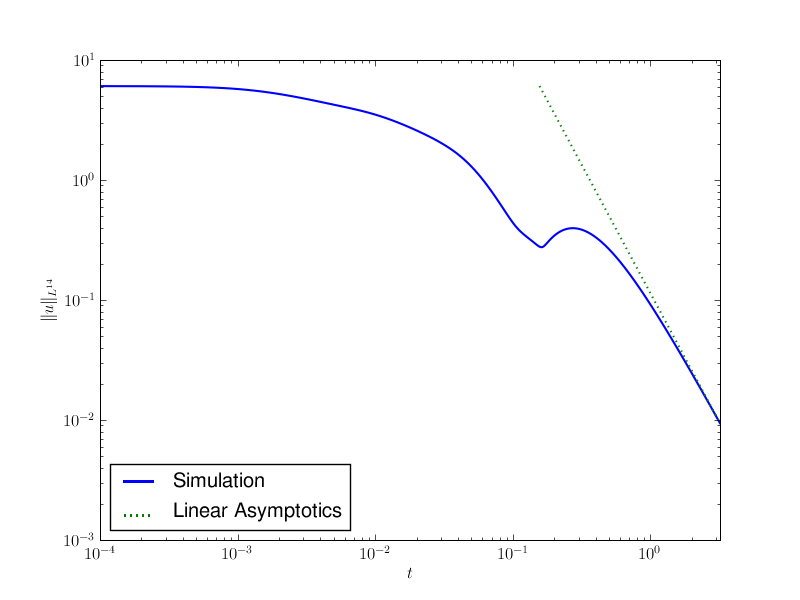}
\\
 &\\
 \mathrm{(e)} & \mathrm{(f)}
 \end{array}$
\caption{Metrics for $u_0 = 10 e^{-r^2}$.  Figures (b), (d), (f) are computed on the domain $[0, 2000]$, with 200000+1 points.}
%(a) and (b) show the saturation of $\dot H^2$.  (a) includes measurements of$\dot B^{2}_{2,\infty}$.   (c) and (d) show the decay of $\norm{u}_{L^6}$.  (e) reveals spreading in the Fourier domain.  (f) plots the evolution of the $L^{14}_x$ norm and the anticipated asymptotic decay as the system scatters.  Figures (a), (c), and (e) are computed on the domain $[0, 100]$ with $40000+1$ points.  
\label{fig:gauss10metrics}
\end{figure}

\begin{figure}
$\begin{array}{cc}
\includegraphics[width=3.5in]{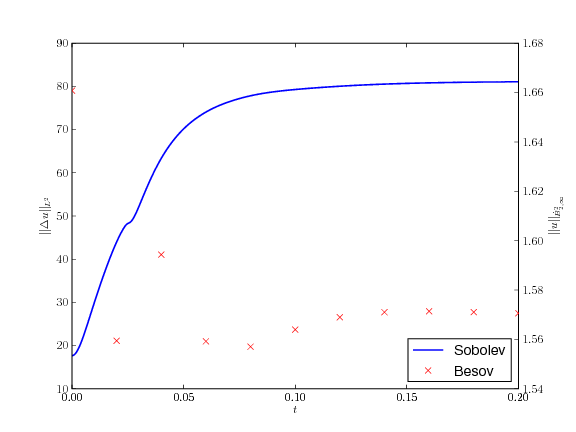}  &
\includegraphics[width=3.5in]{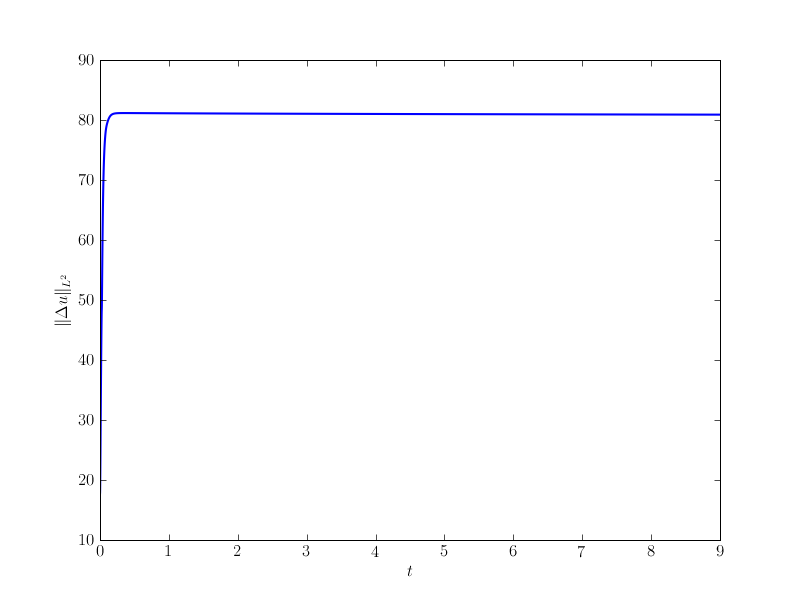}\\
\mathrm{(a)} & \mathrm{(b)}\\
\includegraphics[width=3.5in]{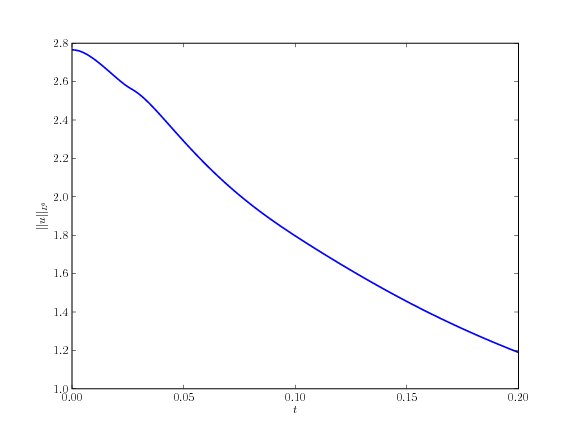}  &
\includegraphics[width=3.5in]{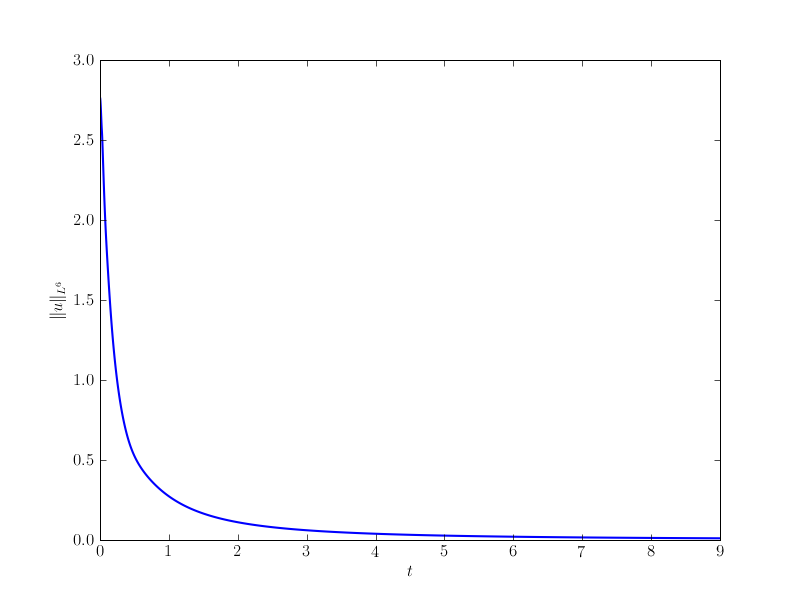}\\
\mathrm{(c)} & \mathrm{(d)}\\
\includegraphics[width=3.5in]{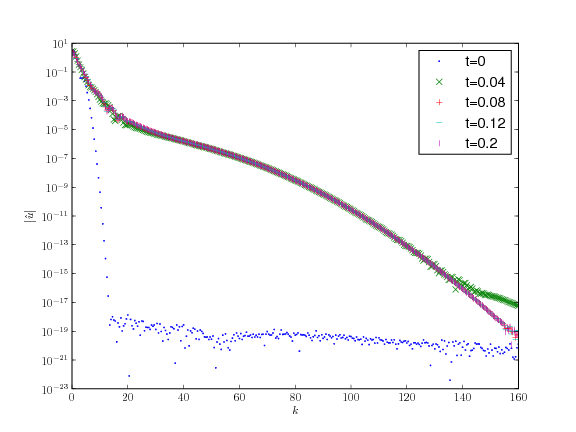}  &  \includegraphics[width=3.5in]{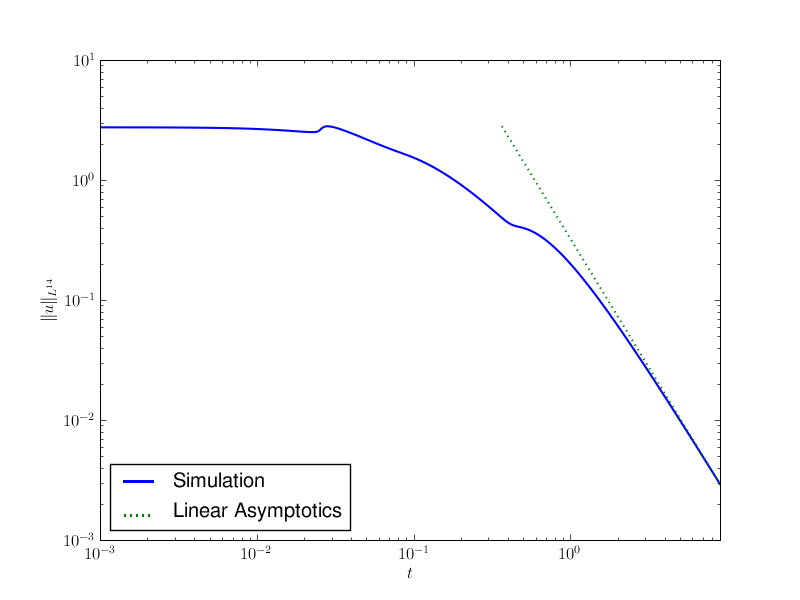}
\\
 &\\
 \mathrm{(e)} & \mathrm{(f)}
 \end{array}$
\caption{Evolution of $u_0 =8 r^2 e^{-r^2}$.  Figures (a), (c), and (e) are computed on the domain $[0, 100]$ with $32000+1$ points.  Figures (b), (d), and (f) are computed on the domain $[0, 2400]$, with 120000+1 points.}
\label{fig:ring7metrics}
\end{figure}

\begin{figure}
$\begin{array}{cc}
\includegraphics[width=3.5in]{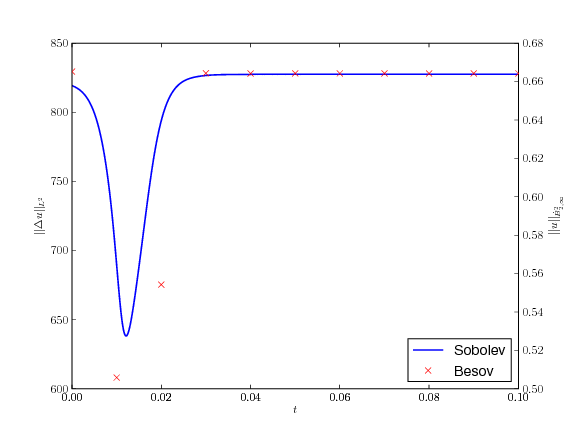}  &
\includegraphics[width=3.5in]{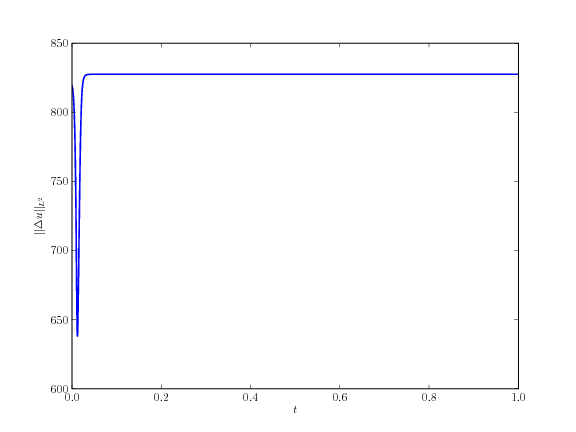}\\
\mathrm{(a)} & \mathrm{(b)}\\
\includegraphics[width=3.5in]{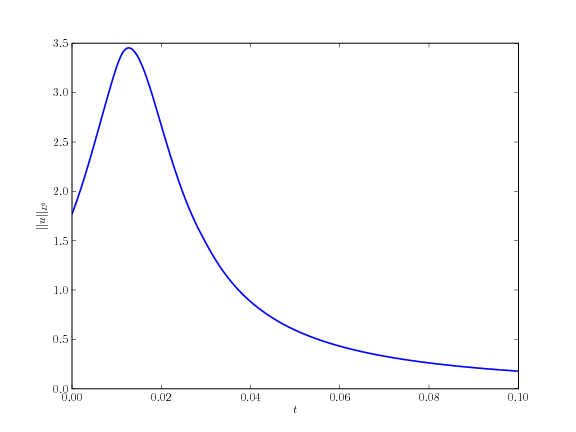}  &
\includegraphics[width=3.5in]{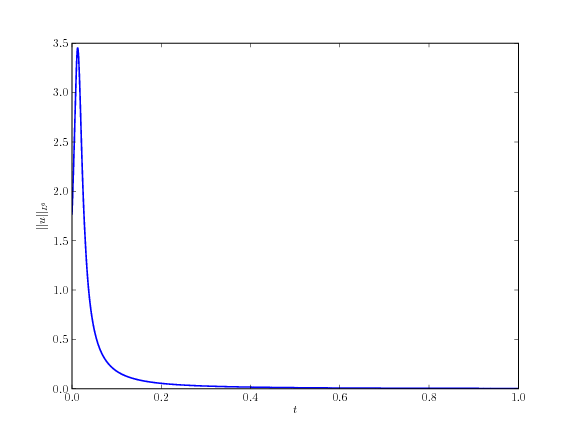}\\
\mathrm{(c)} & \mathrm{(d)}\\
\includegraphics[width=3.5in]{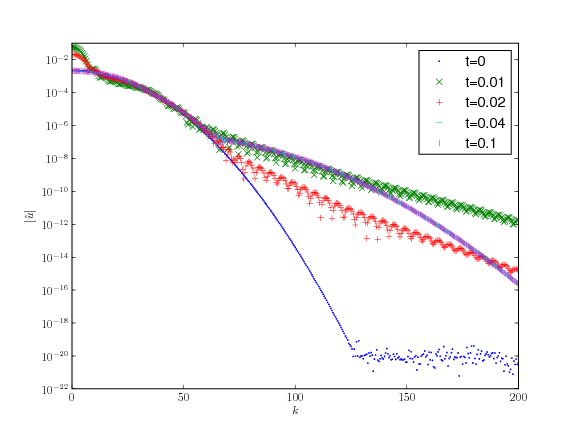}  &  \includegraphics[width=3.5in]{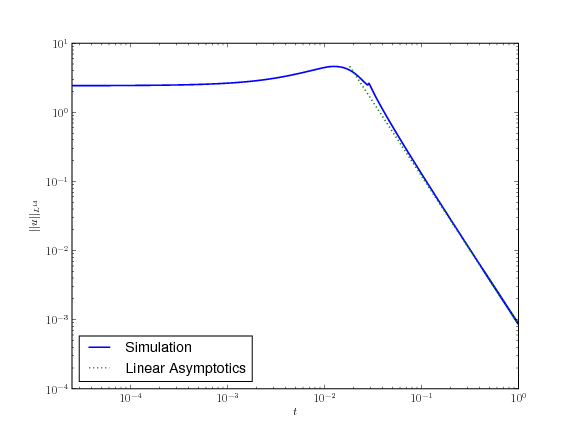}
\\
 &\\
 \mathrm{(e)} & \mathrm{(f)}
 \end{array}$
\caption{Evolution of $u_0 = 4 e^{-10\mathrm{i}r^2} e^{-r^2}$.  Figures (a), (c), and (e) are computed on the domain $[0, 100]$ with $40000+1$ points.  Figures (b), (d), and (f) are computed on the domain $[0, 1000]$, with 200000+1 points.}
\label{fig:ring5metrics}
\end{figure}

Much of the recent analytical progress on NLS used careful treatment of the equation in the Fourier domain.  We examine the Fourier transform of our simulations for hints that might be applied to future analysis.  The transform, plotted at various times in Figures \ref{fig:gauss10metrics} (e), \ref{fig:ring7metrics} (e), and \ref{fig:ring5metrics} (e) shows several features.  There is an initial spreading into high wave numbers, and the support is much broader than the initial condition.  This relaxes, and the limiting state has a smaller support than during the transient period, but still in excess of the initial condition.  The asymptotic constancy of the transform is further evidence that $u$ evolves linearly and scatters.

We can also interpret the Fourier data through the Besov norm \eqref{besov-def}.  Using the method described in Section \ref{sec:besov_numerics}, we approximate the Besov norm $\dot{B}^2_{2,\infty}$ at several times in each simulation.  This data, appearing as $\times$'s in Figures \ref{fig:gauss10metrics} (a),  \ref{fig:ring7metrics} (a), and \ref{fig:ring5metrics} (a), shows several things.  First, the Besov norm $\dot{B}^2_{2,\infty}$ is orders of magnitudes smaller than the scale invariant norm $\dot{H}^2$.  Like the Sobolev norm, there is some transient variation followed by saturation to some asymptotic value.  In the case of the oscillating gaussian data, Figure \ref{fig:ring5metrics} (a), the dynamics of the two norms appear to be phase locked.  Another feature is that while the Sobolev norm can increase (substantially) as it saturates, the Besov norm always decays.  More detailed information for each of the three simulations is given in Tables \ref{table:gauss10besov2}, \ref{table:ring7besov2}, and \ref{table:ring5besov2}.

\begin{table}
\begin{center}
\begin{tabular}{r | l | l }
Time & $\norm{\Delta u}_{L^2}$ & $\norm{u}_{\dot{B}^2_{2,\infty}}$\\
\hline
0.000 & 20.2791& 1.68739\\
0.004 &679.386 &1.12284 \\
0.008 &864.094 & 1.26507\\
0.012 & 1040.59& 1.17468\\
0.016 &1119.92 &1.23148\\
0.020 & 1150.59& 1.26737\\
0.024 & 1163.46& 1.28262\\
0.028 & 1169.47& 1.28782\\
0.032 &1172.59 & 1.28787\\
0.036 & 1174.35& 1.28542\\
0.040 & 1175.41&1.28225
\end{tabular}
\caption{Comparison of norms for the Gaussian data with $R_{\max} = 100$ and $N = 10000$.  Both Besov and Sobev saturate very rapidly.}
\label{table:gauss10besov2}
\end{center}
\end{table}

\begin{table}
\begin{center}
\begin{tabular}{r | l | l }
Time & $\norm{\Delta u}_{L^2}$ & $\norm{u}_{\dot{B}^2_{2,\infty}}$\\
\hline
0.000 & 17.6789& 1.66075\\
0.002 & 43.9913  & 1.55944\\
0.004 &63.3487 &1.59434 \\
0.006 & 74.0549 & 1.55922  \\
0.008 &77.8274 & 1.55706\\
0.010 & 79.2784 & 1.56396  \\
0.012 & 80.0381& 1.56898\\
0.014 & 80.5317  & 1.57103  \\
0.016 & 80.8207 &1.57138\\
0.018 & 80.9811 & 1.57106\\
0.020 & 81.0696& 1.57057
\end{tabular}
\caption{Comparison of norms for the ring data with $R_{\max} = 100$ and $N = 32000$.  Both Besov and Sobev saturate very rapidly.}
\label{table:ring7besov2}
\end{center}
\end{table}

\begin{table}
\begin{center}
\begin{tabular}{r | l | l }
Time & $\norm{\Delta u}_{L^2}$ & $\norm{u}_{\dot{B}^2_{2,\infty}}$\\
\hline
0.00 & 819.277& 0.665221\\
0.01 & 689.723  & 0.505825\\
0.02 &793.98 &0.554248 \\
0.03 & 826.662 & 0.664362  \\
0.04 &827.449 & 0.6642\\
0.05 & 827.483 & 0.664262  \\
0.06 & 827.486& 0.664268\\
0.07 & 827.487  & 0.664269  \\
0.08 & 827.487 &0.664269\\
0.09 & 827.487 & 0.664269\\
0.1 & 827.487& 0.664269
\end{tabular}
\caption{Comparison of norms for the oscillating gaussian data with $R_{\max} = 100$ and $N = 40001$.  Both Besov and Sobev saturate very rapidly.}
\label{table:ring5besov2}
\end{center}
\end{table}

\bibliography{energy_nls}
\newpage

\end{document}